\documentstyle[12pt,amstex,amssymb]{article}
\begin{document}
\baselineskip15pt
\textwidth=12truecm
\textheight=20truecm
\hoffset-.1cm
\voffset-.5cm

\title{Superselection Theory for Subsystems}
\author{
Roberto Conti\\
Dipartimento di Matematica, Universit\`a di Roma ``Tor Vergata,''\\
I-00133 Rome, Italy\\
\\
Sergio Doplicher\\
Dipartimento di Matematica, Universit\`a di Roma ``La Sapienza,''\\
I-00185 Rome, Italy\\
\\
John E.\ Roberts\\
Dipartimento di Matematica, Universit\`a di Roma ``Tor Vergata,''\\
I-00133 Rome, Italy\\
}
\date{}
\maketitle

\begin{abstract}
An inclusion of observable nets satisfying duality induces an inclusion
of canonical field nets. Any Bose net intermediate between the observable
net and the field net and satisfying duality is the fixed--point net of
the field net under a compact group. This compact group is its canonical
gauge group if 
the occurrence of sectors with
infinite\\ statistics can be ruled out for
the observable net and its vacuum Hilbert space is separable.
\end{abstract}

\vfill

\noindent Research supported by MURST, CNR--GNAFA and EU.\eject

\section{Introduction}

   In this paper, we take the view that a physical system is
described by its observable net, a net $\cal{O}\mapsto\frak{A}(\cal{O})$
of von Neumann algebras over double cones in Minkowski space in its
vacuum representation. Our goal is to analyze subsystems. Of course,
physical intuition would lead us to believe that in physically realistic
situations there will be no proper subsystems since any putative subsystem
loses its identity through interaction with the ambient system. A result
in this direction would be interesting as it would allow one to pinpoint
many natural sets of generators. For example, one could, modulo 
technicalities,
claim that the observable net is generated by the energy--momentum tensor
density \cite{Dop},\cite{Con}.
Unfortunately, we have no such result.
Instead, we turn to consider
what might prove to be the exceptional case where there are proper 
subsystems.
For example, if the original system admits symmetries, i.e.\ if there is a
nontrivial
group of local automorphisms of the net leaving the vacuum state invariant,
then the fixed--point net provides an example of a subsystem and one may 
even wonder whether every subsystem arises in this manner.\smallskip

   In fact, the theory of superselection sectors provides a natural
mechanism for giving examples of subsystems. Each observable net $\frak A$
is contained in an associated canonical field net $\frak{F}$ as the
fixed--point net under a compact group $G$ of gauge automorphisms. The
fixed--point net $\frak{B}$ under a proper subgroup $G$, containing the
element changing the sign of Fermi fields, if present, can be treated as
the observable net of some other physical system and $\frak{A}$ will then
be a subsystem of $\frak{B}$. \smallskip

   Here are some of the questions we would like to answer. Can one classify
subsystems? What is the relation between the superselection structure of
a system and that of its subsystems and how are their canonical field nets
related? Earlier partial results on the classification problem
can be found in \cite{LaSc} and \cite{Dav}.\smallskip

   We have only been able to make sensible progress on classifying subsystems
by restricting ourselves to systems satisfying duality in the vacuum
representation. Since there are good reasons for believing that $\frak{A}$
satisfies essential duality, this amounts to replacing $\frak{A}$ by its
dual net, thus ignoring, for example, the possibility of spontaneously broken
gauge symmetries. This partial solution does have the merit of reducing the
classification problem to that of finding all observable nets with a given
dual net.\smallskip

   In Section 2 we study inclusions of observable nets and their functorial
properties. Thus if $\frak{A}\subset\frak{B}$, we get an inclusion of the
corresponding categories of $1$--cocycles, 
$Z^1({\frak A})\subset Z^1({\frak B})$
and hence of the corresponding categories of transportable
morphisms, ${\cal T}_t({\frak A})\subset{\cal T}_t({\frak B})$, and
finally restricting to finite statistics gives
${\cal T}_f({\frak A})\subset{\cal T}_f({\frak B})$. We interpret this
latter
inclusion in terms of the associated homomorphism from the gauge group of
$\frak{B}$ to that of $\frak{A}$.\smallskip

  In Section 3, we study inclusions of field nets giving conditions for the
existence of an associated conditional expectation. Conditional expectations
also play a decisive role in proving that an inclusion of observable algebras
satisfying duality give rise to an inclusion of the corresponding complete
normal field nets.\smallskip

   In Section 4, we study intermediate nets, that is nets contained
between the observable net and its canonical field net, showing that
such nets are the fixed--points of $\frak{F}$ under a closed subgroup $L$
of the field net. After this, we prove a result showing that the sectors
of an intermediate observable net correspond, as one would expect, to
the equivalence classes of irreducible representations of $L$. This
was the principal objective of this paper and is worth comparing with
previous results. There are two known sets of structural hypotheses 
\cite{Dri},
\cite{Re} allowing one to conclude that a Bosonic net has no sectors.
Our result would be a consequence whenever $\frak{F}$ were Bosonic. However,
in the absence of evidence that a Bosonic canonical field net satisfies the
structural hypotheses, these results of \cite{Dri} and \cite{Re}
have been most useful in proving the absence of sectors in examples such
as free field theories. For our result, we need to exclude infinite
statistics for the observable net, a weaker hypothesis with little 
known about its validity. In addition we have to assume that
the vacuum Hilbert space of $\frak{A}$ is separable, as indeed it is, 
in practice.\smallskip

  The paper concludes with an appendix giving results on the harmonic analysis
of the action of compact groups on von Neumann and $C^*$--algebras and on
conditional expectations needed in the course of this paper.

\section{Inclusions of Observable Nets}

  In this section, we will be considering an inclusion 
${\frak A}\subset{\frak B}$
of observable (i.e. local) nets,
with a view to seeing what can be said about the relation between the
corresponding superselection sectors.
Each observable net will be considered as acting
irreducibly on its own vacuum Hilbert space
(denoted in the sequel by ${\cal H}_{\frak A}$, ${\cal H}_{\frak B}$).
Of course, if ${\frak A}\subset{\frak B}$, ${\cal H}_{\frak A}$
is naturally identified with a subspace of ${\cal H}_{\frak B}$.

  However, we start by recalling some well known facts about superselection
structure, cf.\ \cite{H}, \cite{K}, \cite{DHRa}, \cite{DHRb}. This will
allow us to introduce our notation and give a few definitions and useful
results.

Throughout this paper, the superselection sectors for ${\frak A}$
are understood to be
the unitary equivalence classes of irreducible representations $\pi$ satisfying
the Selection Criterion (SC)
$$\pi |_{{\frak A}({\cal O}')} \cong
\pi_0 |_{{\frak A}({\cal O}')}, \ {\cal O} \in {\cal K}$$
with respect to
a reference vacuum representation $\pi_0 \ (=\pi^{\frak A}_0).$
Here ${\cal K}$ denotes as usual the set of double
cones in Minkowski space ordered under inclusion.

The representations satisfying the selection criterion are the objects of
a W${}^*$-category $S(\pi_0)$ whose arrows are
the intertwining operators.

Recall that ${\frak A}^d$, the dual net of ${\frak A}$,
is defined by
${\frak A}^d({\cal O}):={\frak A}({\cal O}')' $,
the commutants being taken on ${\cal H}_{\frak A}$.
If $\frak A$ is a local net,
i.e.\ if ${\frak A}\subset {\frak A}^d$,
we have inclusions
${\frak A} \subset {\frak A}^{dd}:=({\frak A}^d)^d
 \subset {\frak A}^d={\frak A}^{ddd}.$

If ${\frak A}$ satisfies Haag duality, i.e.\ if ${\frak A}={\frak A}^d$,
we find that a $\pi$ as above is equivalent to a representation of
the form $\pi_0 \circ \rho,$ where $\rho$ is an
endomorphism of $\frak A$ with the following properties:
it is localized in a double cone ${\cal O},$
i.e.\ $\rho\upharpoonright_{{\frak A}({\cal O}')}={\rm id},$
and is transportable, i.e.\ (inner) equivalent
to an endomorphism localized in any other double cone. The
category of these endomorphisms and their intertwiners
is equivalent to $S(\pi_0)$ and hence may be equally used to describe
the superselection structure. However, it has the added advantage of being
a tensor $W^*$--category, denoted by ${\cal T}_t$, and reveals the latent
tensor structure of superselection sectors. The objects of ${\cal T}_t$,
the set of localized and transportable endomorphisms of $\frak A$,
will be denoted by $\Delta_t({\frak A})$.\smallskip

   If we wish to relate the superselection sectors of ${\frak A}$ and
${\frak B}$ using endomorphisms
we run into two evident problems.
On one hand, the restriction of an endomorphism $\rho$ of ${\frak B}$
to ${\frak A}$ is not an endomorphism of ${\frak A},$
unless $\rho(\frak A)\subset {\frak A}$, although it could still be
regarded as a representation.
Furthermore, it is not at first sight clear how
a given localized, transportable
endomorphism of $\frak A$ can be extended to a similar endomorphism
of $\frak B.$
However, this problem has a canonical solution
indicated by an alternative approach to superselection sectors
using net cohomology.\smallskip

   In the version in \cite{Ro}, net cohomology is conceived as a cohomology
of partially ordered sets with coefficients in nets over the
partially ordered set $\cal K$.
The formal description of this cohomology may be found in \cite{Ro}
and we restrict ourselves here to a pedestrian account of those
concepts needed in this paper.\smallskip

   A $0$--simplex $a$ is just a double cone ${\cal O}$, i.e.\ an element of
${\cal K}$. A $1$--simplex $b$ is an ordered set $({\cal O},{\cal O}_0,
{\cal O}_1)$ of double cones with  ${\cal O}_0 \cup {\cal O}_1 \subset
{\cal O}$.
Its faces are the $0$--simplices $\partial_ob={\cal O}_0$,
$\partial_1b={\cal O}_1$ and its support is $|b|={\cal O}$.
A $2$--simplex $c$ is an ordered set of double cones
$({\cal O},{\cal O}_0,{\cal O}_1, {\cal O}_2, {\cal O}_{01},{\cal O}_{02},
{\cal O}_{12})$ such that ${\cal O}\supset {\cal O}_0\cup{\cal O}_1
\cup{\cal O}_2$ and such that its faces $\partial_0c=({\cal O}_0,
{\cal O}_{01},{\cal O}_{02})$, $\partial_1c=({\cal O}_1,{\cal O}_{01},
{\cal O}_{12})$ and $\partial_2c=({\cal O}_2,{\cal O}_{02},{\cal O}_{12})$
are $1$--simplices. Its support is $|c|={\cal O}$. The set of
$n$--simplices is denoted by $\Sigma_n$.

\noindent{\bf Definition} {\sl
A $0$--cocycle of ${\frak A}$, a net of $C^*$--algebras over
${\cal K}$, is a map $z: \Sigma_0 \to {\frak A}$ such that
$z(\partial_0 b)=z(\partial_1 b)$, $b \in \Sigma_1,$
and $z(a)\in {\frak A}(a)$, $a \in \Sigma_0$.}\smallskip

Hence the set $Z^0({\frak A})$ of $0$--cocycles is
$\cap_{\cal O}{\frak A}({\cal O}).$

\noindent{\bf Definition} {\sl
A $1$--cocycle of  ${\frak A}$ is a map $z: \Sigma_1 \to
{\cal U}({\frak A})$ such that
$z(\partial_0 c)z(\partial_2 c)=z(\partial_1 c)$, $c \in \Sigma_2$,
and $z(b)\in {\frak A}(|b|),$ $b\in\Sigma_1$.}\smallskip

   The $1$--cocycles are considered as the objects of a
$C^*$--category $Z^1({\frak A})$. An arrow between
$1$--cocycles, $w\in(z,z')$ is a mapping $w: \Sigma_0 \to{\frak A}$
such that $z(b)w(\partial_1 b)=w(\partial_0 b)z'(b),$ $b \in \Sigma_1,$
$w(a) \in {\frak A}(a), a \in \Sigma_0 .$
Note that if $1$ denotes
the trivial $1$-cocycle $1(b)=I, b \in \Sigma_1$, then the elements
of $(1,1)$ are just the $0$--cocycles. Two objects
$z, z'$ of $Z^1({\frak A})$ are cohomologous if $(z,z')$ contains a
unitary arrow and $z$ is a $1$--coboundary
if it is cohomologous to $1.$

 Here is an example of $1$--cocycle illustrating at the same time the
relation with the theory of superselection sectors. Given a
representation $\pi$ of ${\frak A}$ satisfying the selection criterion,
pick for each $a\in\Sigma_0$ a unitary operator $V_a$ such that
$$V_a\pi(A)=\pi_0(A)V_a,
\quad A\in{\frak A}({\cal O}),\,\,a\subset{\cal O}',$$
and set
$$z(b):=V_{\partial_0b}V_{\partial_1b}^*,\quad b\in\Sigma_1,$$
then $z\in Z^1({\frak A}^d)$.\smallskip

  Conversely, given a $1$--cocycle with values in a net ${\frak A}$,
define for $a\in\Sigma_0$,
$$\pi_a(A):=z(b)Az(b)^*,\quad 
\text{provided}\,\,A\in{\frak A}^d({\cal O})\,\,,
b\in\Sigma_1,\,\partial_0b=a,\,\partial_1b\subset{\cal O}'.$$
One checks that $\pi_a$ gives a well defined representation of
${\frak A}^d$ and that $z(b)\in(\pi_{\partial_1b},\pi_{\partial_ob})$.
Thus a $1$--cocycle gives rise to a field $a\mapsto \pi_a$ of
equivalent representations. Furthermore, $\pi_a$ is localized in $a$
in the sense that
$$\pi_a(A)=A,\quad A\in{\frak A}^d({\cal O}),\,\,{\cal O}\subset a'.$$
Details may be found in \cite{Ro}, \S 3.4.6, Theorem 1, Corollary 2,
where it is also proved that $S(\pi_0)$
and $Z^1({\frak A}^d)$ are equivalent $W^*$--categories. It follows
that $S(\pi_0)$ and $S(\pi^{dd}_0)$ are equivalent as $W^*$--categories,
where $\pi^{dd}_0$ denotes the vacuum representation of the double
dual net ${\frak A}^{dd}$.\smallskip

   It should be noted that the above results on superselection
sectors do not require any form of duality or even locality.
But we are not able to define the tensor structure without a further
hypothesis. We see, however, that essential duality,
${\frak A}^d={\frak A}^{dd}$, will suffice for this purpose.\smallskip

   A variant of the above construction relates cocycles and
endomorphisms. It is based on assuming relative duality
$${\frak A}({\cal O})={\frak A}^d({\cal O})\cap{\frak A},
\quad {\cal O}\in{\cal K},$$
a weaker version of the more familiar assumption of duality. This is
the $C^*$--version of duality and is defined without reference to
the vacuum representation. In this context, it is natural to use nets
of $C^*$--algebras. As a consequence of relative duality,
an endomorphism localized in ${\cal O}$ satisfies
$\rho({\frak A}({\cal O}_1))\subset{\frak A}({\cal O}_1)$ whenever
${\cal O}\subset{\cal O}_1$. Furthermore, an intertwiner between
endomorphisms localized in ${\cal O}$ is automatically in
${\frak A}({\cal O})$. Relative duality suffices for a theory
of transportable endomorphisms but to pass from superselection
sectors to transportable endomorphisms we need duality or, at
least, essential duality.\smallskip

   Consider an endomorphism $\rho$ of ${\frak A}$ such that, given
$a\in\Sigma_0$, there is a unitary $\psi(a)\in{\frak A}$ with
$$\psi(a)\rho(A)=A\psi(a)\quad A\in{\frak A}({\cal O}),\,\, {\cal O}\subset
a'.$$
This is the analogue for endomorphisms of the selection criteria for
representations. Our $1$--cocycle
$z(b):=\psi(\partial_ob)\psi(\partial_1b)^*$,
$b\in\Sigma_1$ now takes values in ${\frak A}$. Any such $1$--cocycle
$z$ now defines a field of endomorphisms:
$$\rho_a(A):=z(b)Az(b)^*,\quad 
\text{provided}\,\,A\in{\frak A}({\cal O})\,\,,
b\in\Sigma_1,\,\partial_0b=a,\,\partial_1b\subset{\cal O}'.$$
$\rho_a$ is localized in $a$ in the sense that
$$\rho_a(A)=A,\quad A\in{\frak A}({\cal O}),\,\,{\cal O}\subset a'.$$
Since $z(b)\in(\rho_{\partial_1b},\rho_{\partial_0b})$, we
have a field $a\mapsto\rho_a$ of endomorphisms in ${\cal T}_t$, each of
which is equivalent to the $\rho$ we started from. Note that if we
start with a cocycle of the form
$z(b):=\psi(\partial_ob)\psi(\partial_1b)^*$, as above, then
$\rho_a=\text{Ad}\psi_a\rho$. In particular, if $\rho$ is localized
in $a$ we may take $\psi_a=I$ and hence arrange that $\rho_a=\rho$.
We may regard our construction as leading to an equivalence of tensor
$C^*$--categories between $Z^1({\cal T}_t)$ and ${\cal T}_t$,
cf.\ \cite{Ro}, \S 3.4.7, Theorem 5.\smallskip

   As in the theory of superselection sectors, the tensor $C^*$--category
${\cal T}_t({\frak A})$ admits a canonical permutation symmetry
$\varepsilon$ (in more than two spacetime dimensions).
\medskip

   We can now begin to examine an inclusion ${\frak A}\subset{\frak B}$
of nets satisfying relative duality. Such an inclusion obviously induces
an inclusion functor $Z^1({\frak A}) \to Z^1({\frak B}).$ Thus a
$1$--cocycle in a local net ${\frak A}$ not only gives rise to a field
$a\mapsto\rho_a$ of endomorphisms of ${\frak A}$ but to a field
$a\mapsto\tilde\rho_a$ of endomorphisms of any relatively local net
${\frak B}$ extending the original field. We have seen that any element
of $\Delta_t({\frak A})$ arises as a value of such a field and hence
admits an extension to an element of $\Delta_t({\frak B})$. As the
cocycle is not uniquely determined by the endomorphism, a little
argument is needed to show that the extension is uniquely detemined,
cf.\ Lemma 3 of \S 3.4.7 in \cite{Ro}.\smallskip

\noindent
{\bf Lemma 2.1} {\sl Let $z$ and $z'$ be two $1$--cocycles 
of a net ${\frak A}$
satisfying relative duality and suppose that, for some $a\in\Sigma_0$,
$\rho_a=\rho'_a$. Then, if ${\frak A}\subset{\frak B}$ is an inclusion of
nets and ${\frak A}$ and ${\frak B}$ are relatively local, the
endomorphisms $\tilde\rho_a$ and $\tilde\rho'_a$ of ${\frak B}$
induced by $z$ and $z'$ agree.}\smallskip

\noindent
{\bf Proof.} Let $b\in\Sigma_1$ with $\partial_0b=a$ then
$z(b)^*z'(b)\in(\rho'_{\partial_1b},\rho_{\partial_1b})$ is
an intertwiner of endomorphisms localized in $\partial_1b$.
Since ${\frak A}$ satisfies relative duality,
$z(b)^*z'(b)\in{\frak A}(\partial_1b)$ and the result follows
since ${\frak B}$ and ${\frak A}$ are relatively local.\smallskip

   We now come to the main result of this section.\smallskip

\noindent
{\bf Theorem 2.2} {\sl Let ${\frak A}\subset{\frak B}$ be an inclusion
of nets satisfying relative duality. Then there is an induced
structure preserving inclusion of ${\cal T}_t({\frak A})$ in
${\cal T}_t({\frak B})$ which corresponds to the above extension
on endomorphisms and to the given inclusion on intertwiners.}\smallskip

\noindent
{\bf Proof.} In view of the relation between cocycles and endomorphisms
and Lemma 2.1, the only point which is not yet obvious is that the tensor
structure is preserved by the inclusion. However, if $z$ and $z'$ are
$1$--cocycles in ${\frak A}$ and $a\mapsto\rho_a$ and $a\mapsto\rho'_a$
are the corresponding fields of endomorphisms, then
$$z\otimes z'(b):=z(b)\rho_{\partial_1b}(z'(b))$$
defines a $1$--cocycle over ${\frak A}$ whose associated field
of endomorphisms is $a\mapsto\rho_a\rho'_a$, and the result now
follows.\smallskip

The extension of endomorphisms is also discussed in \cite{BE} under the
name $\alpha$--induction in the context of nets of subfactors \cite{LR}.

   The inclusion of Theorem 2.2 will of course map the unitary
operator $\varepsilon(\rho,\rho')$ in ${\cal T}_t({\frak A})$
onto the corresponding operator for the extended endomorphisms.
Furthermore, as is obvious from the cohomological description,
we shall have
$$(\rho,\rho')_{\frak{A}}=(\rho,\rho')_{\frak{B}}\cap\frak{A},$$
with an obvious notation.\smallskip

   In particular, whenever ${\frak A}$ and ${\frak B}$ satisfy duality,
this result is at the same time a result about superselection structure
and relates the superselection structure of ${\frak A}$ to that of
${\frak B}$. The extension of an endomorphism $\rho$ with finite
statistics, $\rho\in\Delta_f({\frak A})$, will again have finite
statistics and we have an induced tensor $^*$--functor from
${\cal T}_f({\frak A})$ to ${\cal T}_f({\frak B})$. In fact
this also holds in the more general context provided we
understand ${\cal T}_f$ to be the full subcategory of ${\cal T}_t$
having conjugates. Now ${\cal T}_f({\frak A})$ and ${\cal T}_f({\frak B})$
are equivalent to the tensor $W^*$--categories of finite dimensional
continuous unitary representations of compact groups so that
tensor $^*$--functors correspond contravariantly to continuous
homomorphisms between the groups in question\cite{DRf}. In the context of
superselection structure, the compact groups are the gauge groups.
A gauge group appears as the group of automorphisms of a field net
leaving the observable subnet pointwise fixed.\smallskip

   We would like to make this homomorphism explicit and therefore
consider the following situation. We consider a commuting square
of inclusions of nets ${\frak A}_1\subset{\frak A}_2\subset{\frak F}_2$
and ${\frak A}_1\subset{\frak F}_1\subset{\frak F}_2$. The ${\frak A}_i$
are to be considered as observable nets, the
${\frak F}_i$ as field nets, cf. \cite{DR}, Definition 3.1 and
Theorem 3.6. Thus we
suppose ${\frak A}_i$ to have trivial relative commutant in ${\frak F}_i$
and ${\frak F}_i$ to be local relative to ${\frak A}_i$. We consider
the subcategories ${\cal T}_i$ of ${\cal T}_f({\frak A}_i)$ induced
by Hilbert spaces in ${\frak F}_i$. The Hilbert spaces in question are
unique and are supposed to generate
${\frak F}_i$. We let $G_i$ be the group of automorphisms of
${\frak F}_i$ leaving ${\frak A}_i$ pointwise invariant. These
automorphisms leave the Hilbert spaces stable and we suppose that
$G_i$ is a compact group equipped with the topology of pointwise
norm convergence on these Hilbert spaces. Finally, we suppose that
each irreducible representation of $G_i$ is realized on some
such Hilbert space and that ${\frak F}_i^{G_i}={\frak A}_i$. These
last conditions ensure that the inclusion ${\frak A}_i\subset{\frak F}_i$
realizes ${\cal T}_i$ in a canonical way as a dual of $G_i$.
Without wishing to get involved in further technicalities, we might
say that the essence of these conditions is that the net ${\frak F}_i$
is a crossed product of the net ${\frak A}_i$ by the action
of a group dual, where these terms are to be understood as adaptions
to nets of von Neumann algebras of the corresponding concepts in
\cite{DRe}.\smallskip

\noindent
{\bf Theorem 2.3} {\sl Under the above conditions on a commuting
square of inclusions, ${\frak F}_1$ is stable under the action of
$G_2$ and the restriction of $G_2$ to ${\frak F}_1$ defines a
homomorphism $h$ from $G_2$ to $G_1$. The inclusion
${\frak A}_1\subset{\frak A}_2$ induces an inclusion functor
from ${\cal T}_1$ to ${\cal T}_2$ and this inclusion functor
is precisely that induced by $h$. If $N$ and $K$ denote the
kernel and image of $h$, respectively, then
$${\frak F}_1\vee{\frak A}_2={\frak F}_2^N,$$
$${\frak F}_1\cap{\frak A}_2={\frak F}_1^K.$$
}

\noindent
{\bf Proof.} Note first that a Hilbert space $H(\rho)$ in ${\frak F}_1$
inducing an object $\rho$ of ${\cal T}_1$ must, when considered
as a Hilbert space in ${\frak F}_2$, induce the canonical extension
of $\rho$ to an object of ${\cal T}_f({\frak A}_2)$ by relative
locality. This canonical extension is thus an object of
${\cal T}_2$ so that we do have an induced inclusion functor
from ${\cal T}_1$ to ${\cal T}_2$. It also follows that
$H(\rho)$ is stable under the action of $G_2$. But
such Hilbert spaces generate ${\frak F}_1$ so ${\frak F}_1$
is stable under the action of $G_2$. The restriction of an
element of $G_2$ to ${\frak F}_1$ defines an automorphism of ${\frak F}_1$
leaving ${\frak A}_1$ pointwise invariant and is therefore
an element of $G_1$. Thus restriction defines the required homomorphism
$h$. Since the representation of $G_2$ on $H(\rho)$ arises by
composing that of $G_1$ with $h$, $h$ induces the above inclusion
functor. Now $N$, being the kernel of $h$, obviously acts trivially
on ${\frak F}_1$ and ${\frak A}_2$. Now ${\frak F}_2^N$ is generated
by the Hilbert spaces $H(\rho)$ in ${\frak F}_2$
inducing objects of ${\cal T}_2$ and carrying irreducible
representations of $G_2$ that are trivial in restriction to $N$.
Regarding these as representations of $K$ and inducing up to a
representation of $G_1$,
bearing in mind that every irreducible
representation of $G_1$ is realized within ${\frak F}_1$, we conclude
that there is an isometry in ${\frak A}_2$ mapping $H(\rho)$ into
${\frak F}_1$. Thus ${\frak F}_2^N$ is generated by ${\frak F}_1$
and ${\frak A}_2$. Next, note that the $K$--invariant part of a
Hilbert space of ${\frak F}_1$ inducing an object of ${\cal T}_1$
is $G_2$--invariant and hence lies in ${\frak A}_2$.  These
Hilbert spaces generate ${\frak F}_1^K$ and, as any element of
${\frak F}_1\cap{\frak A}_2$ is $K$--invariant, we have
${\frak F}_1\cap{\frak A}_2={\frak F}_1^K$, completing the
proof.\smallskip

\section{Inclusions of Field Nets}

  In the last section, we have treated inclusions of observable nets. However,
observable nets are frequently defined by starting with a net ${\frak F}$
of fields with Bose--Fermi commutation relations. From a mathematical point 
of view, these are simply the ${\Bbb Z}_2$--graded version of an
observable net.
Hence to have a basic formalism which is sufficiently flexible, we need to
consider inclusions of ${\Bbb Z}_2$--graded nets.\smallskip

   We define a
(concrete) ${\Bbb Z}_2$--graded net ${\frak F}$ to be a net of von Neumann
algebras over ${\cal K}$, represented on its (vacuum) Hilbert
space ${\cal H}_{{\frak F}}$, together with an involutive unitary operator
$k$ inducing a net automorphism $\alpha_k$ of ${\frak F}$. The even (Bose) 
part ${\frak F}_+$  of ${\frak F}$ is the fixed--point net under
$\alpha_k$, the odd (Fermi) part ${\frak F}_-$ changes sign under
$\alpha_k$.
The twisted net ${\frak F}^t$ is defined as ${\frak F}_++ik{\frak F}_-$
and is, in an obvious way, itself a ${\Bbb Z}_2$--graded net.
The ${\Bbb Z}_2$--graded or twisted dual net of ${\frak F}$ is defined by
$${\frak F}^d({\cal O}):=\cap_{{\cal O}_1\subset {\cal O}'}
{\frak F}^t({\cal O}_1)'.$$
It is understood to act on the same Hilbert space with the same unitary $k$.
${\frak F}$ satisfies twisted duality if it coincides
with its twisted dual net ${\frak F}^d$.
${\frak F}$ is said to have Bose--Fermi commutation relations if
$${\frak F}({\cal O}_1)\subset{\frak F}^t({\cal O}_2)',\quad {\cal O}_1
\subset {\cal O}_2',$$
or, equivalently, if ${\frak F}\subset{\frak F}^d$. If, in addition,
${\frak F}$
is irreducibly represented on ${\cal H}_{{\frak F}}$,
we refer to the triple ${\frak F}, k,{\cal H}_{{\frak F}}$ as being a field
net.\smallskip

   By an inclusion of ${\Bbb Z}_2$--graded nets we mean compatible 
(normal) inclusions ${\frak B}({\cal O})\subset{\frak F}({\cal O})$ of von
Neumann algebras together with an inclusion of Hilbert spaces 
${\cal H}_{{\frak B}}\subset{\cal H}_{{\frak F}}$ compatible 
with the inclusion of nets and such that $k_{{\frak B}}$ is the
restriction of $k_{{\frak F}}$ to ${\cal H}_{{\frak B}}$.
We further require that ${\cal H}_{{\frak B}}$ be cyclic and separating
for each ${\frak F}({\cal O})$.\smallskip

   Typically, an observable net may be defined from a
field net acted on by a compact group $G$ of net automorphisms with
$\alpha_k\in G$
by taking ${\frak A}$ to be the fixed--point net under $G$. Under these
circumstances, ${\frak A}$ satisfies duality if ${\frak F}$ satisfies
twisted duality (except in one space dimension) and there is a normal
conditional expectation $m$ of nets from ${\frak F}$ onto ${\frak A}$
obtained by averaging over the group. However, it has been known since
the beginnings of the theory of superselection sectors that the existence
of such a normal conditional expectation follows simply from the hypothesis
that ${\frak A}$ satisfies duality, without any reference to a compact group
$G$. We present here some related results.\smallskip

Let ${\frak F}, k,{\cal H}_{{\frak F}}$ be a ${\Bbb Z}_2$--graded net and
let $E$ be a projection on ${\cal H}_{{\frak F}}$, commuting with $k$ and
cyclic and separating for ${\frak F}$, i.e.\ for each ${\frak F}({\cal O})$.
Let ${\frak F}^E({\cal O}):={\frak F}({\cal O})\cap \{E\}'$ and let
${\frak F}_E$ and $k_E$ denote the restriction of ${\frak F}^E$ and $k$
to the subspace $E{\cal H}_{{\frak F}}$. Then the triple
${\frak F}_E,k_E,E{\cal H}_{{\frak F}}$ is itself a ${\Bbb Z}_2$--graded net.
If we started with a field net, we would only get a field net if we knew
that ${\frak F}_E$ acts irreducibly on $E{\cal H}_{{\frak F}}$.\smallskip

   We now ask whether ${\frak F}$ admits a conditional expectation $m$ of
nets such that
$$m(F)E=EFE,\quad F\in{\frak F}.$$
In this case, $m$ would project onto the subnet ${\frak F}^E$ and be locally
normal, see Lemma A.7 of the Appendix. In particular, in the case of a field
net ${\frak F}_E$ would act irreducibly on $E{\cal H}_{{\frak F}}$.\smallskip

\noindent{\bf Lemma 3.1} {\sl If ${\frak F}$ is a field net and ${\frak F}_E$
satisfies twisted duality, then there is a conditional expectation of
${\frak F}$ such that
$$m(F)E=EFE,\quad F\in{\frak F}.$$}\smallskip

\noindent{\bf Proof} By Corollary A.8b of the Appendix, we must show that
$$[EFE,EF'E]=0,\quad F\in{\frak F}({\cal O}),\quad F'\in{\frak F}({\cal O})'.$$
Now if ${\cal O}_1\subset{\cal O}'$ and $B\in{\frak F}^t({\cal O}_1)_E$
then $[EFE,B]=0$. Hence $E{\frak F}({\cal O})E\upharpoonright E{\cal
H}_{{\frak F}}
\subset({\frak F}_E)^d({\cal O})={\frak F}_E({\cal O})$. Thus there is a
$G\in{\frak F}({\cal O})$ with $GE=EG=EFE$, and $[GE,EF'E]=0$, as
required.\smallskip

  To have a more systematic approach, we begin by proving an analogue of
Lemma A.7 of the Appendix for ${\Bbb Z}_2$--graded nets.\smallskip

\noindent{\bf Lemma 3.2} {\sl Let ${\frak F},k$ be a ${\Bbb Z}_2$--graded net
on a Hilbert space ${\cal H}$. Let $E$ be a $k$--invariant projection cyclic
for ${\frak F}$ and ${\frak F}^d$. Let ${}^E{\frak F}$ be the net defined by:
$${}^E{\frak F}({\cal O}):=\{F\in{\frak F}({\cal O}):EFE\in(E{\frak F}^dE)^d
({\cal O})\}.$$
Then ${}^E{\frak F}\supset{\frak F}^E$ and is weak--operator
closed. This makes ${}^E{\frak F}$ into a ${\frak F}^E$--bimodule. Given
$F\in{\frak F}({\cal O})$, there is a $m(F)\in{\frak F}^{ddE}({\cal O})$
such that
$$m(F)E=EFE$$
if and only if $F\in{}^E{\frak F}({\cal O})$.}\smallskip

\noindent{\bf Proof} If $F\in{}^E{\frak F}({\cal O})$, then
$$EFE\in (E{\frak F}^{dt}E)({\cal O}_1)',\quad {\cal O}_1\subset{\cal O}'.$$
Pick $G\in{\frak F}^{dt}({\cal O}_1)$, then
$$EF^*EG^*GEFE\leq \|EFE\|^2EG^*GE,$$
hence there exists 
$m(F)\in{\frak F}^{dt}({\cal O}_1)'$, such that $m(F)E=EFE$.
Since $E$ is separating for each ${\frak F}^{dt}({\cal O}_1)'$ 
and ${\cal O}'$
is path--connected, $m(F)$ is independent of the choice of ${\cal O}_1
\subset{\cal O}'$. Hence $m(F)\in\cap_{{\cal O}_1\subset{\cal O}'}
{\frak F}^{dt}({\cal O}_1)'={\frak F}^{dd}({\cal O}_1)$. If
$F\in{\frak F}({\cal O})$ and there is an $m(F)\in{\frak F}^{ddE}({\cal O})$,
such that $m(F)E=EFE$ then $F\in{}^E{\frak F}({\cal O})$. The
remaining assertions are evident.\smallskip

\noindent{\bf Remark} To have a closer analogy with Lemma A.7 of the Appendix,
we should require that ${\frak F}={\frak F}^{dd}$. Since, after all,
we require ${\cal M}={\cal M}''$ in the Appendix.\smallskip

\noindent{\bf Corollary 3.3} {\sl Let ${\frak F}={\frak F}^{dd},k,{\cal H}$
be a ${\Bbb Z}_2$--graded net and $E$ a projection cyclic for ${\frak F}$
and ${\frak F}^d$. Then the following conditions are equivalent.
\begin{description}
\item{a)} There is a conditional expectation $m$ on ${\frak F}$ such that
$m(F)E=EFE$, $F\in{\frak F}$.
\item{a$'$)} There is a conditional expectation $m^d$ on ${\frak F}^d$
such that $m^d(F^d)E=EF^dE$, $F^d\in{\frak F}^d$.
\item{b)} $E{\frak F}E\subset (E{\frak F}^dE)^d$.
\item{b$'$)} $E{\frak F}^dE\subset (E{\frak F}E)^d$.
\item{c)} $({}_E{\frak F})^d={}_E({\frak F}^d)$.
\item{c$'$)} $(_E({\frak F}^d))^d={}_E{\frak F}$.
\end{description}
Here ${}_E{\frak F}$, for example, denotes the restriction of $E{\frak F}E$ to
$E{\cal H}$.}\smallskip

\noindent{\bf Proof} Suppose b) holds, then ${}^E{\frak F}={\frak F}$
and, by Lemma 3.2, $m$ becomes a conditional expectation onto ${\frak F}_E$,
since it is idempotent and of norm 1, giving a). Similarly, b$'$) implies
a$'$). Taking duals, we see
that b) and b$'$) are equivalent. It is clear that
a) implies b) and that a$'$) implies b$'$). Now
$$(_E{\frak F})^d({\cal O})=\cap_{{\cal O}_1\subset{\cal O}'}
(_E{\frak F}^t)({\cal O}_1)'=\cap_{{\cal O}_1\subset{\cal O}'}(E{\frak F}^t
({\cal O}_1)E\upharpoonright E{\cal H})',$$
so if a) holds then by Corollary A.8c of the Appendix,
$$(_E{\frak F})^d({\cal O})=\cap_{{\cal O}_1\subset{\cal O}'}
(E{\frak F}^t({\cal O}_1)'E)\upharpoonright E{\cal H}=
E\cap_{{\cal O}\subset{\cal O}'}{\frak F}^t({\cal O}_1)'E\upharpoonright
E{\cal H},$$
where we have used the fact that $E$ is separating for each
${\frak F}^t({\cal O}_1)'$ and that ${\cal O}'$ is path--connected.
Thus a) implies c). The implication a$'$) implies c$'$) follows by
exchanging the role
of ${\frak F}$ and ${\frak F}^d$. Now, trivially, c) implies b$'$)
and, again, c$'$) implies b) follows.\smallskip

   Of course, a direct application of Corollary A.8 of the Appendix
shows that the above conditions are also equivalent to
$$({}_E{\frak F})({\cal O})'={}_E({\frak F}({\cal O})'),\quad {\cal
O}\in{\cal K}.$$
However, in view of the superficial similarities with c), it is worth 
stressing that equivalence depends on being in more than two spacetime
dimensions.\smallskip

   What becomes clear from the above discussion is that the problem of
studying the subsystems of a given system can be divided up in a natural
way. We can begin with the simple class of subsystems characterized by
cyclic projections $E$ and the existence of a conditional expectation as above.
Let us call such subsystems {\it full} since they are the largest subsystems
on their Hilbert spaces and are uniquely determined by their Hilbert
spaces. If ${\frak A}$ is a full subsystem of ${\frak B}$ and ${\frak B}$ is
itself a full subsystem of ${\frak F}$, then ${\frak A}$ is a full
subsystem of ${\frak F}$. We see from Lemma 3.1 that a subsystem satisfying
twisted duality is full. Furthermore, if ${\frak F}$ satisfies twisted
duality, then, by Corollary 3.3,  a subsystem is full if and only if it
satisfies twisted duality. A second step might then be to analyze
subsystems having the same Hilbert space.\smallskip

   In the following result, we give an analogue of Corollary A.9 of the
Appendix and look at full subsystems from the point of view of the
subsystem.\smallskip

\noindent {\bf Lemma 3.4} {\sl Let ${\frak B}\subset{\frak F}$ be an inclusion
of ${\Bbb Z}_2$--graded nets and $E$ the associated projection from ${\cal
H}_{{\frak F}}$
to ${\cal H}_{{\frak B}}$, then the following conditions are equivalent.
\begin{description}
\item{a)} There is a (necessarily unique, injective and ${\Bbb Z}_2$--graded)
net morphism $\nu:{\frak B}^d\to{\frak F}^d$ such that
$$\nu(B)\Phi=B\Phi,\quad B\in{\frak B}^d,\quad \Phi\in{\cal H}_{{\frak B}}.$$
\item{b)} ${\frak B}^d={\frak F}^d{}_E$.
\end{description}
If the conditions are fulfilled, there is a unique normal conditional
expectation $m$ of ${\frak F}^d$ onto $\nu({\frak B}^d)$ such that
$$m(F)E=EFE,\quad F\in {\frak F}^d.$$}\smallskip

\noindent{\bf Proof} $\nu$ is obviously unique, hence ${\Bbb Z}_2$--graded,
since ${\cal H}_{{\frak B}}$ is cyclic for each ${\frak F}^t({\cal O})$,
hence separating for each ${\frak F}^d({\cal O})$.
Given a), we note that $E\nu(B)E=\nu(B)E$ and replacing $B$ by $B^*$, we see
that $\nu(B)\in{\frak F}^d{}^E$. Hence $B\in{\frak F}^d{}_E$,
yielding b). 
Conversely, if b) is satisfied, given $B\in{\frak B}^d({\cal O})$,
there is an $F\in{\frak F}^d({\cal O})$ with $FE=EF$ and $F\Phi=B\Phi$,
$\Phi\in{\cal H}_{{\cal B}}$. Hence,
we may pick $\nu(B)=F$ to give a map $\nu:{\frak B}^d\to{\frak F}^d$ and it
follows
from uniqueness that $\nu$ is a net morphism. Now suppose the conditions are
satisfied and that $F\in{\frak F}^d({\cal O})$ and $B\in{\frak B}^t({\cal
O}_1)$
with ${\cal O}_1\subset{\cal O}'$. Then
$$EFEB=EFBE=EBFE=BEFE.$$
Hence the restriction of $EFE$ to ${\cal H}_{{\frak B}}$ lies in
${\frak B}^d({\cal O})={\frak F}^d{}_E({\cal O})$ by b). The result now
follows by Lemma A.7 of the Appendix.\smallskip

\noindent{\bf Remarks} For an inclusion of field nets, ${\frak
B}\subset{\frak F}^d{}_E
\subset{\frak B}^d$, b) is trivially fulfilled if ${\frak B}$ satisfies
twisted duality. Now suppose that ${\frak B}$ satisfies twisted
duality for wedges then
$${\frak B}^d({\cal O}) = \cap_{{\cal W}\supset{\cal O}}{\frak R}({\cal W}),$$
where ${\frak R}({\cal W})$ denotes the von Neumann algebra associated with
the wedge ${\cal W}$. Now given spacelike double cones, ${\cal O}$ and
${\cal O}_1$, there is a wedge ${\cal W}$ such that ${\cal O}\subset{\cal W}
\subset{\cal O}_1'$. Hence
$${\frak B}^d({\cal O})\subset{\frak R}({\cal W})\subset{\frak F}^t({\cal
O}_1)'{}_E,$$
and, taking the intersection over ${\cal O}_1$, we see that b) is again
satisfied. If ${\frak B}$ satisfies essential twisted duality, i.e.\
if ${\frak B}^d={\frak B}^{dd}$, then we cannot conclude from the
above that b) is satisfied since we do not know that we have an inclusion
${\frak B}^{dd}\subset{\frak F}^{dd}$. If ${\frak B}={\frak B}^{dd}$,
${\frak F}={\frak F}^{dd}$ and $E$ is also cyclic for each
${\frak F}^d(\cal{O})$
in Lemma 3.4, then we may deduce from Corollary 3.3 that, under the equivalent
conditions of Lemma 3.4, ${\frak B}={\frak F}_E$.\smallskip

   We now consider an inclusion $\frak A\subset\frak B$ of nets of local
von Neumann algebras over double cones each satisfying duality in their
respective Hilbert spaces $\cal H_{\frak A}$ and $\cal H_{\frak B}$. Let
$E$ denote the projection of $\cal H_{\frak B}$ onto $\cal H_{\frak A}$.
Then, as follows e.g.\ from Lemma 3.4, there is a conditional expectation of
nets of von Neumann algebras $m$ of $\frak B$ onto $\frak A$ such that
$$EBE=m(B)E,\quad B\in\frak B.$$
Furthermore, the intertwiners spaces between transportable localized
morphisms of $\frak{A}$ and their extensions to $\frak{B}$ are related by
$$m(\rho,\rho')_{\frak{B}}=(\rho,\rho')_{\frak{A}},$$\smallskip
see the remarks following Theorem 2.2.

We now introduce the canonical field net $\frak{F}$ of $\cal B$ and let
$m_{\frak B}$ be the associated conditional expectation from $\frak F$ onto 
$\frak B$. We recall\cite{DR} that the canonical field net is defined for
observable nets satisfying duality and Property B.
Let $\frak{E}$ denote the net of $C^*$--algebras generated by
the Hilbert spaces in $\frak F$ implementing the transportable localized
morphisms of $\frak A$. Then by Lemma A.2, the restriction of $m\circ
m_{\cal B}$ to $\frak{E}$ is the unique conditional expectation 
$n$ 
onto the subnet $\frak A$.
By \cite{DRe}, this shows that $\frak{E}$ is the $C^*$--cross product of
$\frak A$
by the action of $\cal T_{\frak A}$.\smallskip

  Now let $\alpha$ denote the canonical action of the gauge group $G$ of
$\frak A$ on the net $\frak{E}$. Then we have
$$\int\alpha(F)\,d\mu(g)=n(F),\quad F\in\frak{E},$$
where $\mu$ denotes Haar measure on $G$. Let $\cal H_{\frak F}$ denote the
canonical Hilbert space of $\frak F$ and $\cal H$ the Hilbert subspace
generated by $\cal H_{\frak A}$ and $\frak{E}$. Since\\
$\omega\circ m\circ m_{\frak B}=\omega$,
$$ \omega\circ n=\omega$$
where $\omega$ is a state defined by a vector of $\cal H_{\frak A}$.
It follows that states of $\frak{E}$ defined by vectors
in $\cal H_{\frak A}$ are gauge invariant. Hence
$$(\Phi,F\Psi)=(\Phi,n(F)\Psi),\quad F\in\frak{E},\quad \Phi,\Psi\in\cal
H_{\frak A}.$$
Given $F_i\in\frak{E}$ and
$\Phi_i\in\cal H_{\frak A}$, $i=1,2\dots,n$ define
$$U_g\sum_i F_i\Phi_i:=\sum_i\alpha_g(F_i)\Phi_i,$$
$$\|\sum_i\alpha_g(F_i)\Phi_i\|^2=\sum_{i,j}(\Phi_i,\alpha_g(F_i^*F_j)\Phi_j)
=\sum_{i,j}(\Phi_i,F_i^*F_j\Phi_j).$$
Thus we get a unitary action of $G$ on $\cal H$.\smallskip

   We next remark that $\cal H_{\frak A}$ is the space of $G$--invariant
vectors in $\cal H$. In fact, if $U_g\Phi=\Phi$, $g\in G$ and $\Phi$ is
orthogonal to $\cal H_{\frak A}$. then
$$(\Phi,F\Psi)=(\Phi,\alpha_g(F)\Psi)=(\Phi,n(F)\Psi)=0,\quad
\Psi\in\cal H_{\frak A},\quad F\in\frak{E}.$$
Thus $\Phi=0$.\smallskip

   We can now check easily, that our data consisting of a representation
of $\frak A$ on $\cal H$ restricting to the vacuum representation on
$\cal H_{\frak A}$, a unitary action of $G$ on $\cal H$ and a homomorphism
$\rho\mapsto H_\rho$ from the semigroup of objects of $\cal T_{\frak A}$
has all the properties needed to generate the canonical field net of
$\frak A$ \cite{DR}, p.66. Since $\cal H_{\frak A}$ is cyclic for $\frak
E(\cal O)'$
hence separating for $\frak E(\cal O)''$, the canonical field net is
canonically isomorphic to $\cal O\mapsto\frak E(\cal O)$. Thus we have shown
the following result.\medskip

\noindent{\bf Theorem 3.5} {\sl Let $\frak A\subset\frak B$ be an inclusion
of nets of observable algebras satisfying duality and Property B, 
then there is a canonical inclusion of the corresponding canonical field
nets.}

\section{Sector Structure of Intermediate Nets}

   In this section, we consider an inclusion of nets ${\frak A}\subset
{\frak B}$
and examine in more detail the relation between the sectors of ${\frak A}$ and
those of ${\frak B}$.
As we have little to say in general, we restrict our attention to the case
that ${\frak A}\subset{\frak B}\subset \frak {F(A)}$. We first show that
under these circumstances, $\frak{B}$ is the fixed--point net of $\frak{F(A)}$
under a closed subgroup of the gauge group $G$ of $\frak{A}$. To this end,
we denote by $H_\rho$ the Hilbert space in $\frak{F}:=\frak{F(A)}$
inducing $\rho\in \Delta_f$. Set $K_\rho:= H_\rho\cap\frak{B}$. Then
$K_\rho$ is a Hilbert space in $\frak{B}$. We claim
\begin{description}
\item{a)} $K_\rho K_\sigma\subset K_{\rho\sigma}$,
\item{b)} $TK_\rho\subset K_\sigma$, if $T\in(\rho,\sigma)$,
\item{c)} $K_{\bar\rho}=JK_\rho$, where $J$ is an antiunitary from $H_\rho$
to $H_{\bar\rho}$ intertwining the actions of the gauge group.
\end{description}
Indeed a) is obvious whilst b) follows from the fact that $T\in\frak{A}$.
Finally, c) follows from the fact that we may define such an antiunitary
$J$ by $J\psi=\psi^*\bar R$, with $\bar R\in(\iota,\rho\bar\rho)$ 
as in the definition of conjugate endomorphisms, cf. 
Theorem 3.3 of
\cite[II]{DHRb}, or a standard solution of the conjugate equations,
cf. \cite{LoRo}.
It follows\cite{R} that there is a unique closed subgroup $L$ of the gauge
group
$G$ such that each $K_\rho$ is precisely the fixed-points of the action of
$L$ on $H_\rho$. \medskip

   We now make use of the fact that when $\frak{B}$ satisfies duality,
there is a locally normal conditional expectation $m$ from $\frak{F}$ onto
$\frak{B}$. Let $\psi,\psi'\in H_\rho$, then
$$m(\psi)B=\rho(B)m(\psi),\quad B\in\frak{B}.$$
Hence
$$\psi^*m(\psi')B=\psi^*\rho(B)m(\psi')=B\psi^*m(\psi'),$$
and since ${\frak B}'\cap{\frak F}({\frak A})={\Bbb C} I$,
$\psi^*m(\psi')\in{\Bbb C} I$ and
$m(\psi')\in H_\rho$. Since $\frak{F}$ is generated as a net of linear
spaces closed in say the $s$--topology by the elements of the Hilbert spaces
$H_\rho$, $\frak{B}$ is generated in the same way by $K_\rho$. Thus $\frak{B}$
is the fixed--point net under the action of $L$. Thus we have proved the
following result.\smallskip

\noindent {\bf Theorem 4.1} {\sl Let $\frak{F}$ be the canonical field net
of the
observable net $\frak{A}$ and $\frak{B}$ an intermediate net,
$\frak{A}\subset\frak{B}\subset\frak{F}$, satisfying duality, then there 
is a closed subgroup $L$ of the gauge group $G$ of $\frak{A}$ such that
${\frak B}={\frak F}^L$.}\smallskip

Related results in the context of inclusions of von Neumann algebras
can be found in \cite{ILP}.

\noindent {\bf Lemma 4.2} {\sl The following are equivalent:

\begin{description}
\item{a)} $L$ is a normal subgroup of $G$,
\item{b)} $\alpha_g({\frak B}) \subset {\frak B}$, $g\in G$,
\item{c)} ${\frak B}$ is generated by Hilbert spaces inducing
endomorphisms in $\Delta_f({\frak A})$.
\end{description}

}\smallskip

\noindent
{\bf Proof.} a) $\Rightarrow$ b) is obvious.
Hilbert spaces inducing endomorphisms in $\Delta_f({\frak A})$
are $G$-invariant so c) $\Rightarrow$ b).
If b) holds then given $g \in G$ and $k \in L$, $B \in {\frak B}$
$$\alpha_{gkg^{-1}}(B)=\alpha_g\alpha_k\alpha_{g^{-1}}(B)
=\alpha_g\alpha_{g^{-1}}(B)=B$$
since $\alpha_{g^{-1}}(B)\in {\frak B}$. Thus $\alpha_{gkg^{-1}}$
is an automorphism of ${\frak F}$ leaving ${\frak B}$
pointwise fixed. Thus $gkg^{-1} \in L$ and $L$ is a normal subgroup,
giving b) $\Rightarrow$ a).
Suppose a) then consider the set of Hilbert spaces in
${\frak B}$ inducing endomorphisms in $\Delta_f({\frak A})$.
These must be $L$-invariant and each thus carry a canonical
representation of $G/L$ and ${\frak B}^{G/L}={\frak F}^G
={\frak A}$. We know that ${\frak B}$ is generated by
$H \cap {\frak B}$ where $H$ is a Hilbert space inducing an
element of $\Delta_f({\frak A})$. This Hilbert space
may not have support $I$ but it is an invariant subspace for
the action of $G$ and is hence in the algebra generated
by Hilbert spaces above. Obviously, c) implies b),
completing the proof.
\smallskip

We next discuss a situation where two members of an inclusion of observable 
nets ${\frak A} \subset {\frak B}$ have coinciding canonical field nets.
We start with a net $\frak B$ and suppose that its canonical field net
$\frak F$ has a compact gauge group $K$ of internal symmetries with
$K \supset G,$ where $G$ is the gauge group of $\frak B$ and then define
$\frak A$ to be the fixed--point net ${\frak F}^K$.

We recall that if $K$ is spontaneously broken then ${\frak A}$ does not
satisfy duality\cite{R}. Its dual net ${\frak A}^d$ is the fixed-point net of
${\frak F}$ under the closed subgroup of unbroken symmetries and does
satisfy duality. Furthermore, ${\frak A}$ and ${\frak A}^d$ have the same
superselection structure. Hence in line with our strategy of
considering only nets satisfying duality, we may restrict ourselves
to the case that $K$ is unbroken. We recall that,
if ${\frak F}$ has the split property, then the group $K^{max}$
of all unitaries leaving $\Omega$ invariant and inducing net automorphisms
of $\frak F$ is automatically compact in the strong operator 
topology\cite{DLb}.

In the above situation $\{ {\frak F},{\frak A},K, {\cal H}_{\frak A} \}$
is a field system with gauge symmetry for ${\frak A}.$
Furthermore, $\rho \in \Delta_f({\frak A})$
is induced by a finite-dimensional Hilbert space $H$ in ${\frak F}$
since this is true of its extension to an element of $\Delta_f({\frak B}).$
But this means that every sector of ${\frak A}$ is realized on the vacuum
Hilbert space of ${\frak F}$
so that ${\frak F}$ is the canonical field net of
${\frak A}$ and $K$ is the gauge group. We have thus proved the following
result\smallskip

\noindent
{\bf Proposition 4.3} {\sl Let $\frak{B}$ be an observable net with
canonical field net $\frak{F}$ and gauge group $G$. Suppose $\frak{F}$ has
an unbroken compact group $K$ of internal symmetries. Then the fixed--point
net ${\frak F}^K$ has $\frak{F}$ as canonical field net.}\smallskip

  Finally, we consider the sector structure of an intermediate observable
net $\frak{A}\subset\frak{B}\subset\frak{F}(\frak{A})$ satisfying duality.
As we know from Theorem 4.1, 
$\frak{B}$ is the fixed--points of $\frak{F}(\frak{A})$ under the
action of a closed subgroup $L$ of the gauge group $G$. We shall suppose
that the vacuum Hilbert space of $\frak{A}$ is separable, that Property B
of Borchers holds for ${\frak A}^d$ and that 
{\it each representation of $\frak{A}$
satisfying the selection criterion is a direct sum of irreducibles with 
finite statistics}.\smallskip

   We now pick a representation $\hat\pi$ of $\frak{B}$ satisfying the
selection criterion for $\frak{B}$. To analyse this representation, we
choose an associated 1--cocycle $z$ as in \S 2. Since $\frak{B}$ satisfies
duality, $z(b)\in{\frak B}(|b|)\subset{\frak F}(|b|)$. If we consider
$z$ as a $1$--cocycle of $\frak{F}$, it can be used, as discussed in
\S 2, to define representations $\tilde\pi_a$ of $\frak{F}$ where
$$\tilde\pi_a(F):=z(b)Fz(b)^*,\quad F\in{\frak F}({\cal O}),\,b\in\Sigma_1,\,
\partial_0b=a,\,\partial_1b\subset{\cal O}'.$$
Note that $z(b)$ is a Bosonic operator in $\frak{F}$. 
Restricting $\tilde\pi_a$
first to $\frak{B}$ and then to the vacuum Hilbert space of $\frak{B}$
gives the representations $\hat\pi_a$ associated with $z$ considered as
a cocycle of $\frak{B}$. Thus $\hat\pi_a$ is equivalent to $\hat\pi$.

   We now let $\pi$ denote the restriction of some fixed $\tilde\pi_a$
to $\frak{A}$. Now if the vacuum Hilbert space of $\frak{F}$ is non--separable,
then $\pi$ cannot satisfy the selection criterion as its restriction to
each $\frak{A}(\cal{O}')$ is equivalent to a direct sum of uncountably
many copies of the identity representation of $\frak{A}(\cal{O}')$. However,
we shall see that $\pi$ is just a direct sum of representations satisfying
the selection criterion. It suffices to show that any cyclic subrepresentation
satisfies the selection criterion. Such a cyclic representation is, like $\pi$,
locally normal and hence acts on a separable Hilbert space as a consequence
of the following well known result.\smallskip

\noindent
{\bf Lemma 4.4} {\sl Let $\frak{A}$ be an observable net acting on a separable
vacuum Hilbert space and $\omega$ be a locally normal state. Then the GNS
representation $\pi_\omega$ of $\frak{A}$ is separable.}\smallskip

   A proof may be found for example in \S 5.2 of \cite{BW}.\smallskip

\noindent
{\bf Lemma 4.5} {\sl Every cyclic subrepresentation of $\pi$ satisfies the
selection criterion.}\smallskip

\noindent
{\bf Proof} We turn the equivalence of representations
in restriction to $\frak{A}(\cal{O'})$ into a question of the equivalence of
two projections $E_0$ and $F_0$ in the representation
of $\frak{A}(\cal{O'})$ obtained by restricting the vacuum representation
$\hat\pi_0$ of $\frak{F}$ to $\frak{A}(\cal{O}')$. $E_0$ is the projection 
onto the subspace given by the vacuum sector of $\frak{A}$. $F_0$ is
determined as follows. To be able to exploit the Borchers property, we
choose a double cone ${\cal O}_0$ with ${\cal O}^-_0 \subset {\cal O}$,
and a unitary $U$ such that
$$U\pi(A) = \hat\pi_0(A)U, \ A \in {\frak A}({\cal O}'_0),$$
and set $F_0:=UFU^*$, where $F$ corresponds to the
(cyclic) subrepresentation of $\pi$,
$F \in \pi({\frak A})',$ with separable range.

Let
$\sigma, \hat\sigma_0$ and  $\tau, \hat\tau_0$ denote
the restrictions of
$\pi,\hat\pi_0$ to $\frak{A}({\cal O}')$
and ${\frak A}({\cal O}'_0)$, respectively.
Then $U \in (\tau,\hat\tau_0)\subset(\sigma,\hat\sigma_0)$
and $F_0\in(\hat\tau_0,\hat\tau_0)\subset(\hat\sigma_0,\hat\sigma_0)$.
Since $\hat\pi_0$ is, in restriction to $\frak{A}$, a direct sum of
representations satisfying the
selection criterion, $E_0$ has central support $I$ in both
$(\hat\tau_0,\hat\tau_0)$ and $(\hat\sigma_0,\hat\sigma_0)$.

Thus there are projections
$e_0$ and $f_0$ with $e_0 \prec E_0,$ $f_0 \prec F_0$
and $e_0 \simeq f_0 $ in $(\hat\tau_0,\hat\tau_0)$.
Moreover, by Property B for ${\frak A}^d$,
$e_0 \simeq E_0 $ in $(\hat\sigma_0,\hat\sigma_0)$
Thus $E_0$ is equivalent to the subprojection $f_0$ of $F_0$
in $(\hat\sigma_0,\hat\sigma_0)$.
Since $F_0$ is separable and $E_0$
has infinite multiplicity by Property B, we have
$$E_0 \prec F_0 \prec \infty E_0 \simeq E_0.$$
Thus $E_0$ and $F_0$ are equivalent, completing the proof.\smallskip

\noindent
{\bf Corollary 4.6} {\sl $\pi$ is normal on $\frak{A}$, $\tilde\pi_a$ is
normal on $\frak{F}$ and $\hat\pi$ is normal on $\frak{B}$, where
the term normal refers to the vacuum representation of $\frak{F}$.}\smallskip

\noindent
{\bf Proof.} The first statement follows at once from Lemma 4.5, the
second by invoking Theorem A.6 of the Appendix and the third is obvious
since, as we have seen, $\hat\pi$ is equivalent to a subrepresentation
of the restriction $\tilde\pi_a$ to $\frak{B}$.\smallskip

   Now any normal representation of $\frak{B}$ is just a direct sum
of subrepresentations of the defining representation so we have proved
the following result.\smallskip

\noindent
{\bf Theorem 4.7} {\sl Let $\frak{A}$ be an observable net on a separable
Hilbert space whose dual net satisfies Property B and suppose that every
representation of $\frak{A}$ satisfying the selection criterion is a
direct sum of irreducible representations with finite statistics. 
Let ${\frak B}$ be an intermediate observable net satisfying duality,
\ i.e. \ ${\frak A} \subset {\frak B} \subset {\frak F}({\frak A})$
and $L$ the associated compact group as in Theorem 3.1.
Then every representation of $\frak{B}$ satisfying the selection criterion
is a direct sum of sectors with finite statistics and these are labelled
by the equivalence classes of irreducible representations of $L$.}\smallskip

  As a particular case of this, we note that when $\frak{F}$ contains
Fermi elements, then the Bose part of $\frak{F}$ is the fixed--point algebra
of $\frak{F}$ under ${\Bbb Z}_2$ and has precisely two sectors. 
In the case where $\frak A$ has only a finite number of superselection
sectors, the above result is already known, cf. \cite{C},\cite{Mu}. 
Theorem 4.7 has an immediate corollary.\smallskip

\noindent
{\bf Corollary 4.8} {\sl Under the hypothesis of Theorem 4.7, the field nets
of $\frak{A}$ and $\frak{B}$ coincide, 
$\frak{F}(\frak{A})=\frak{F}(\frak{B})$.}
\smallskip

\section{Appendix}

In this appendix we collect together various results  needed
in the course of this paper. They have in common that they do
not involve the net structure but typically the harmonic analysis
of the action of compact groups on von Neumann algebras and
$C^*$--algebras and conditional expectations. The results are
looked at in terms of the structure of the category of
finite-dimensional continuous, unitary representations of the
group rather than the group itself. Consequently, the results
transcend group theory. This degree of generality is not needed
in this paper.\smallskip

\noindent
{\bf Lemma A.1} {\sl Let $m$ be a conditional expectation from the
$^*$--algebra
$\cal B$ onto the $^*$--subalgebra $\cal A$ and $H$ a Hilbert space in
$\cal B$ such that $m(H)\subset H$, then $m$ restricted to $H$ is the
orthogonal projection onto the closed subspace ${\cal A}\cap H$.}\smallskip

\noindent{\bf Proof} If $\psi,\psi'\in H$ then $\psi^*m(\psi)$ is a scalar.
Thus $\psi^*m(\psi')=m(\psi^*m(\psi'))=m(\psi)^*m(\psi')$ so
$m(\psi)^*\psi'=\psi^*m(\psi')$. Hence $m$ restricted to $H$ is selfadjoint
and as it is anyway involutive, it is the orthogonal projection onto
$m(H)={\cal A}\cap H$, as required.\smallskip

   There are some obvious corollaries of this result. Suppose that $\cal B$
is generated by $\cal A$ and a collection $\cal H$ of Hilbert space in $B$
then there is at most one conditional expectation $m$ of $\cal B$ onto
$\cal A$ such that $m(H)\subset H$ for each $H\in\cal H$. If we suppose
that $\cal B$ is a $C^*$--algebra then it suffices if $\cal A$ and
$\cal H$ generate $\cal B$ as a $C^*$--algebra. If $\cal B$ is a
von Neumann algebra and $m$ is normal then it suffices if $\cal A$
and $\cal H$ generate $\cal B$ as a von Neumann algebra. These results
apply in particular to the case where $\cal B$ is the cross product of
$\cal A$ by the action of a dual object of a compact group.
Note, too, that the hypothesis $m(H)\subset H$ is redundant if the
canonical endomorphism of $H$ maps $\cal A$ into itself and if
${\cal A}'\cap {\cal B}=\Bbb C$. Thus 'minimal' or perhaps better
irreducible cross products have a unique mean.\smallskip

\noindent
{\bf Lemma A.2} {\sl Let $\cal A\subset\cal B\subset\cal F$ be inclusions of
$C^*$--algebras and $m_{\cal B}$ a conditional expectation of $\cal F$ onto
$\cal B$. Let $\cal H$ denote a category of Hilbert spaces  in
$\cal F$ each normalizing $\cal B$ and $\cal A$ and such that
$m_{\cal B}(H)\subset H$ for each object $H$ of $\cal H$.
Let $m$ a conditional
expectation of $\cal B$ onto $\cal A$. Suppose that, 
whenever $H$ is an object
of $\cal H$ and $\sigma_H$ the corresponding endomorphism, then
$$\psi\in{\cal A},\quad \psi A=\sigma_H(A)\psi,\quad A\in{\cal A},$$
implies $\psi\in H$. Let $\cal E$ denote the
$C^*$--subalgebra of $\cal F$ generated by $\cal A$ and the objects $H$ of
$\cal H$ then $m \circ m_{\cal B}$
restricted to
$\cal E$  is the unique conditional expectation $n$ of $\cal E$ onto
$\cal A$ with $n(H)\subset H$ for all objects $H$ of $\cal H$.}\smallskip

\noindent
{\bf Proof} The uniqueness of $n$ holds since $\cal E$ is generated by
$\cal A$ and the objects of $\cal H$ and since $n(H)\subset H$ for each
such object $H$. Now taking $n$ to be the restriction of 
$m \circ m_{\cal B}$
to $\cal E$, $n$ is trivially a conditional expectation onto $\cal A$. If
$\psi\in H$, then
$$n(\psi)A=n(\psi A)=\sigma_H(A)n(\psi),\quad A\in\cal A,$$
since $H$ normalizes $\cal A$. Hence $n(\psi)\in H$ by hypothesis, completing
the proof.\smallskip

   By a partition of the identity on a Hilbert space $\cal {H}$
we mean a set $E_i$, $i\in I$ of (self-adjoint) projections with
sum the identity operator. Each element $X\in\cal{B}(\cal{H})$
can then be written $X=\sum_{i,j}E_iXE_j$ with convergence in
say the $s$--topology. The set of elements for which this sum
is finite forms a $^*$--subalgebra ${\cal B}({\cal H})_I$ of
${\cal B}({\cal H})$ which is a direct sum of the subspaces
$E_i{\cal B}({\cal H})E_j$. We let $s_f$ denote the topology
on the $^*$--subalgebra which is the direct
sum of the $s$--topologies on these subspaces.\smallskip

\noindent{\bf Lemma A.3} {\sl Let $E_i$, $i\in I$ be a partition
of the unit on a Hilbert space $\cal{H}$ and $\pi$ a representation
of ${\cal B}({\cal H})_I$ on $\cal{H}_\pi$, continuous in the
$s_f$--topology when $\cal{B}(\cal{H}_\pi)$ is given the $s$--topology.
Then if $\pi(E_i)$, $i\in I$, is a partition of the identity, $\pi$
extends uniquely to an $s$--continuous representation of
$\cal{B}(\cal{H})$}.\smallskip

\noindent{\bf Proof} If $\pi$ extends to an $s$--continuous representation,
again denoted by $\pi$, we must have
$$\pi(X)=\sum_{i,j}\pi(E_iXE_j),\quad X\in\cal{B}(\cal{H}),$$
so any extension is unique. On the other hand, this expression for
$\pi(X)$ is obviously defined on the dense subspace spanned by the
subspaces $\pi(E_i)\cal{H}_\pi$, $i\in I$. Hence, it suffices to show that
$\pi(X)$ is bounded there.
Let $J$ be a finite subset of $I$ and $E_J:=\sum_{j\in J}
E_j$. Then the von Neumann algebra $E_J{\cal B}({\cal H})E_J$ is a
$^*$--subalgebra of ${\cal B}({\cal H})_I$ so that
$$\|\pi(E_JXE_J)\|\leq \|E_JXE_J\|\leq \|X\|,\quad X\in\cal{B}(\cal{H})$$
and $\pi(X)$ is bounded. Computing matrix elements from the dense subspace,
we see that we have a representation of $\cal{B}(\cal{H})$. To see that it
is normal, it suffices to show that its restriction to the compact operators
is non--degenerate. However, its restriction to the compact operators on
each $E_i\cal{H}$ is non-degenerate on $\pi(E_i)\cal{H}_\pi$. But $\pi(E_i)$
is a partition of the identity, so the result follows.\smallskip

\noindent {\bf Remark} Another way of looking at the above result is
that $\cal{B}(\cal{H})$ is the inductive limit of the von Neumann algebras
${\cal B}(E_J{\cal H})$ as $J$ runs over the set of finite subsets of $I$,
ordered under inclusion. The inductive limit is here understood in the
category of von Neumann algebras with normal, but not necessarily
unit--preserving $^*$--homomorphisms.\smallskip

   We now consider a von Neumann algebra $\cal{M}$ and a faithful, normal
conditional expectation $m$ onto a von Neumann subalgebra $\cal{A}$. Consider
$\cal{M}$ as a left $\cal{A}$--module with the $\cal{A}$--valued scalar
product $m(XY^*)$ derived from $m$.\smallskip

\noindent {\bf Lemma A.4} {\sl A representation $\pi$ of $\cal{M}$,
$s$--continuous in restriction to $\cal{A}$ is also $s$--continuous
in restriction to any submodule $\cal{N}$ of finite rank.}\smallskip

\noindent{\bf Proof} When $\cal{N}$ has finite rank, we can find a
finite orthonormal basis $\psi_i$ using the Gram--Schmidt
orthogonalization process. Thus for each $X\in\cal{N}$, we have
$$X=\sum_im(X\psi_i^*)\psi_i.$$
Suppose $X_n\to X$ in the $s$--topology
on $\cal{N}$. Then $\pi(m(X_n\psi_i^*))\to \pi(m(X\psi_i^*))$
and hence $\pi(X_n)\to \pi(X)$ as required.\smallskip

   We will need some variant of this result where $\cal{M}$ is just
a $C^*$--algebra. We could assume that $\cal{N}$ has a finite orthonormal
basis or say assume that it is a finite--rank projective module where
the coefficients can be chosen continuous in the $s$--topology.\smallskip

   To make a bridge between Lemmas A.3 and A.4, we need another structure
related to the notion of hypergroup. We consider a set $\Sigma$ and a
mapping $(\sigma,\tau)\mapsto \sigma\otimes\tau$ from $\Sigma\times\Sigma$
into the set of finite subsets of $\Sigma$. We suppose further that
$\Sigma$ is equipped with an involution (conjugation) 
$\sigma\mapsto\bar\sigma$
with the property that $\rho\in\sigma\otimes\tau$ if and only if
$\tau\in\bar\sigma\otimes\rho$ and if and only if
$\sigma\in\rho\otimes\bar\tau$.
Furthermore there is a distinguished element $\iota\in\Sigma$ such that
$\iota\otimes\sigma$ and $\sigma\otimes\iota$ both consist of the single
point $\sigma$ for each $\sigma\in\Sigma$.\smallskip

   If $\Sigma$ is as above then a $C^*$--algebra $\cal{B}$ will be said to be
$\Sigma$--graded if there are norm--closed linear subspaces $\cal{B}_\sigma$,
$\sigma\in\Sigma$, spanning $\cal{B}$ such that
$\cal{B}_\sigma^*=\cal{B}_{\bar\sigma}$
and if $\cal{B}_\sigma\cal{B}_\tau\subset\cal{B}_{\sigma\otimes\tau}$.
Here $\cal{B}_{\sigma\otimes\tau}$ denotes the norm--closed subspace spanned
by the $\cal{B}_\rho$ as $\rho$ runs over the elements of $\sigma\otimes\tau$.
Note that $\cal{B}_\iota$ is a $C^*$--subalgebra of $\cal{B}$ and that
each $\cal{B}_\sigma$ is a $\cal{B}_\iota$--bimodule.\smallskip

   A representation $\pi$ of a $\Sigma$--graded $C^*$--algebra $\cal{B}$ is a
representation of $\cal{B}$ on a Hilbert space $\cal{H}$ which is a direct
sum of closed linear subspaces $\cal{H}_\sigma$ such that
$$\pi(\cal{B}_\sigma)\cal{H}_\tau\subset\cal{H}_{\sigma\otimes\tau},$$
where $\cal{H}_{\sigma\otimes\tau}$ is defined in the obvious manner.\smallskip

\noindent{\bf Lemma A.5} {\sl Let $\pi$ be a representation of a
$\Sigma$--graded
$C^*$--algebra and $E_\sigma$ the projection on $\cal{H}_\sigma$ then
$$E_\sigma \pi({\cal B})E_\tau\subset
E_\sigma\pi({\cal B}_{\sigma\otimes\bar\tau})E_\tau.$$}

\noindent{\bf Proof}
$\pi(\cal{B}_\rho)\cal{H}_\tau\subset\cal{H}_{\rho\otimes\tau}$.
Thus $E_\sigma\pi({\cal B}_\rho) E_\tau=0$ unless
$\sigma\in\rho\otimes\tau$,
i.e.\ unless $\rho\in\sigma\otimes\bar\tau$.\smallskip

   The obvious example of the above structure is to consider a compact
group $G$
acting on a $C^*$--algebra $\cal{B}$ and to take $\Sigma$ to be the
set of equivalence classes of irreducible, continuous unitary representations
of $G$. We now set ${\cal B}_\sigma:=m_\sigma({\cal B})$, where
$$m_\sigma(B):= \int_G\alpha_g(B)\overline{\chi_\sigma(g)},\quad B\in\cal{B},$$
and $\chi_\sigma$ denotes the normalized trace of $\sigma$.\smallskip

  In the same way, if $(\pi,U)$ is a covariant representation of
$\{\cal{B},\alpha\}$, we get a representation of the $\Sigma$--graded
$C^*$--algebra $\cal{B}$ by using
$$E_\sigma:=\int_G\overline{\chi_\sigma(g)}U(g)$$
to define the closed linear subspace $\cal{H}_\sigma$.\smallskip

   We now put the above results together in the form of a theorem needed in
the body of the text.\smallskip

\noindent{\bf Theorem A.6} {\sl Given a $C^*$--algebra ${\cal B}$ acting
irreducibly
on a Hilbert space ${\cal H}$ and a continuous unitary representation $U$ of
a compact group $G$ inducing an action $\alpha:G\to{\text Aut}(\cal{B})$
on $\cal{B}$ with full Hilbert spectrum, then every representation
of ${\cal B}$ normal on ${\cal B}^G$ is normal on ${\cal B}$}.\smallskip

\noindent{\bf Proof} Let ${\cal A}$ denote the fixed point algebra and
let $E_\sigma$ be as above. Since $\cal{B}$ is irreducible, $E_\sigma$ is
in the weak closure of $\cal{A}$, so that extending $\pi$ to this weak closure
by normality, we have a partition $\pi(E_\sigma)$ of the unit in the
representation space of $\pi$.
 Then by Lemma A.5 above, $E_\sigma{\cal B}E_\tau$
is finite--dimensional as a left ${\cal A}$--module. Since the action has
full Hilbert spectrum, i.e.\ every irreducible representation of $G$ is
realized on some Hilbert space in $\cal{B}$, the argument of Lemma A.4
applies and shows that a
representation $\pi$ of ${\cal B}$ normal on ${\cal A}$ is normal on
each $E_\sigma{\cal B}E_\tau$. The result now follows from Lemma A.3.\smallskip

   We come now to a result on the existence of normal conditional
expectations, beginning with a simple lemma of interest in its own
right.\smallskip

\noindent {\bf Lemma A.7} {\sl Let ${\cal M}$ be a von Neumann algebra on a
Hilbert space ${\cal H}$ and $E$ a cyclic and separating projection for
${\cal M}$.
Let ${\cal M}^E:={\cal M}\cap\{E\}'$ and $^E{\cal M}:=\{M\in{\cal M}:
EME\in(E{\cal M}'E)'\}$, then $^E{\cal M}$ is a weak-operator closed
${\cal M}^E$--bimodule containing ${\cal M}^E$ as a subbimodule. Given
$M\in{\cal M}$ there is a $\mu( M)\in{\cal M}^E$ such that
$$\mu(M)E=EME$$
if and only if $M\in{}^E{\cal M}$.}\smallskip

\noindent {\bf Proof} Given $M\in{}^E{\cal M}$ and $M'\in\cal{M}'$,
then
$$EM^*EM^{'*}M'EME\leq ||EME\|^2\,EM^{'*}M'E,$$
since $EME$ and $EM^{'*}M'E$ commute. Thus $E$ being cyclic for $\cal{M}'$,
there exists a unique bounded operator $\mu(M)$ such that
$$\mu(M)M'E=M'EME.$$
Obviously, 
$\mu(M)\in\cal{M}$ and a computation shows that $\mu(M^*)=\mu(M)^*$.
Setting $M'=I$, it now follows that $\mu(M)$ commutes with $E$. On the other
hand, if $\mu(M)E=EME$ for some $M\in{\cal M}$ then $M\in{}^E{\cal M}$.
The remaining assertions are evident.\smallskip

   Specializing to the case that $^E{\cal M}={\cal M}$ gives the following
result.\smallskip

\noindent {\bf Corollary A.8} {\sl Let $\cal{M}$ be a von Neumann algebra
on a Hilbert space
$\cal{H}$ and $E$ a cyclic and separating projection for $\cal{M}$. Then
the following conditions are equivalent.
\begin{description}
\item{a)} There is a conditional expectation $\mu$ on $\cal{M}$ such that
$\mu(M)E=EME,\quad M\in\cal{M}$.
\item{a')} There is a conditional expectation $\mu'$ on $\cal{M}'$ such
that $\mu'(M')E=EM'E,\quad M'\in\cal{M}'$.
\item{b)} $[E{\cal M}E,E{\cal M}'E]=0$.
\item{c)} $({}_E({\cal M}'))'={}_E{\cal M}$.
\end{description}
Here ${}_E\cal{M}$, for example, denotes the restriction of $E{\cal M}E$
to $E{\cal H}$.
The conditional expectations $\mu$ and $\mu'$ are automatically
normal.}\smallskip

\noindent {\bf Proof} Suppose b) holds then $^E{\cal M}={\cal M}$ and by
Lemma A.7, $\mu$ becomes a normal conditional expectation onto ${\cal M}^E$
since it is idempotent and of norm $1$.
We have therefore deduced a) and by symmetry a$'$). Now suppose a) holds,
then $\mu(\cal{M})$ is just ${\cal M}^E$ and ${\cal M}={}^E{\cal M}$, proving
b). Furthermore, its restriction to $E\cal{H}$ is
$\mu({\cal M})_E$. Thus $({}_E{\cal M})'=\mu({\cal M})'_E$. Since
$\mu({\cal M})=
{\cal M}\cap (E)'$, elements of the form $M'_1+M'_2EM'_3$ with
$M'_i\in\cal{M}'$
form an  $s$--dense $^*$--subalgebra in its commutant and restricting this
to $E\cal{H}$, we have proved c). Trivially, c) implies b), so the conditions
of the corollary are equivalent. \smallskip

\noindent {\bf Remark} If $\sigma$ is an (inner) automorphism of $\cal{B(H)}$,
the above conditions are satisfied by $\sigma(\cal{M})$ and $\sigma(E)$ and the
corresponding conditional expectation is $\sigma\mu\sigma^{-1}$. In particular,
if $\sigma\cal{M}=\cal{M}$ and $\sigma(E)=E$, then
$\mu\sigma=\sigma\mu$.\smallskip

\noindent {\bf Corollary A.9} {\sl Let ${\cal N}\subset{\cal M}$ be an
inclusion of 
von Neumann algebras on Hilbert spaces ${\cal K}$ and ${\cal H}$,
respectively. Let $E$, the projection from ${\cal H}$ onto ${\cal K}$, be 
cyclic and separating for ${\cal M}$, 
then the following conditions are equivalent.
\begin{description}
\item{a)} There is a (necessarily unique and injective) morphism
$\nu:{\cal N}'\to{\cal M}'$ such that
$$\nu(N')\Phi=N'\Phi,\quad N'\in{\cal N}',\quad \Phi\in{\cal K}.$$
\item{b)} ${\cal N}'={\cal M}'{}_E$.
\item{c)} There is a conditional expectation $m$ of ${\cal M}$ onto
${\cal N}$ such that
$$m(M)E=EME,\quad M\in{\cal M}.$$
Here ${\cal M}'_E$ denotes the restriction of ${\cal M}^{'E}$ to
$E\cal{H}=\cal{K}$.
\end{description}}

\noindent{\bf Proof} $\nu$ is obviously unique since ${\cal K}$ is
cyclic for each ${\cal M}$. Given a), we note that $E\nu(N')E=\nu(N')E$
and replacing $N'$ by $N^{'*}$, we see that $\nu(N')\in{\cal M}^{'E}$.
Hence $N'\in{\cal M}'_E$, yielding b). Conversely, if b) is satisfied,
given $N'\in{\cal N}'$, there is an $M'\in{\cal M}'$ with $M'E=EM'$ and
$M'\Phi=N'\Phi$. Hence, we may pick $\nu(N')=M'$ to give a map
$\nu:{\cal N}'\to{\cal M}'$ and it follows from uniqueness that $\nu$
is a morphism. Now if $M\in{\cal M}$ and $N'\in{\cal N}'$, then by a),
$N'$ commutes with the restriction of $EME$ to ${\cal K}$. Thus
$EME\in{\cal N}\subset{\cal M}^E$ and c) follows from Lemma A.7. 
Conversely, if c) holds then ${\cal N}={\cal M}_E$ and b) follows by
calculating commutants.
\smallskip

\noindent{\bf Remark} Of course, when the conditions of Corollary A.9
are satisfied, there is also a conditional expection $m'$ of ${\cal M}'$
onto $\nu({\cal N}')$ such that
$$m'(M')E=EM'E,\quad M'\in{\cal M}'.$$
This follows from Corollary A.8.

\end{document}